# Connecting the Dots: Discovering the "Shape" of Data

By Michelle Feng, Abigail Hickok, Yacoub H. Kureh, Mason A. Porter, Chad M. Topaz

## Abstract

Scientists use a mathematical subject called *topology* to study the shapes of objects. An important part of topology is counting the numbers of pieces and holes in objects, and people use this information to group objects into different types. For example, a doughnut has the same number of holes and the same number of pieces as a teacup with one handle, but it is different from a ball. In studies that resemble activities like "connect the dots", scientists use ideas from topology to study the shape of data. Data can take many possible forms: a picture made of dots, a large collection of numbers from a scientific experiment, or something else. The approach in these studies is called *topological data analysis,* and it has been used to study the branching structures of veins in leaves, how people vote in elections, flight patterns in models of bird flocking, and more. Scientists can take data on the way veins branch on leaves and use topological data analysis to divide the leaves into different groups and discover patterns that may otherwise be hard to find.

## What is a shape?

Shapes are so fundamental to our existence that our brains start to notice them when we are perhaps as young as four or five months. But what exactly do we mean by *shape*? We're used to describing common shapes like lines, circles, and cubes, but what about more complicated objects, like a dragon or a Pokémon or a human being?

*Topology* is a branch of mathematics that studies the shapes of things [3]. To understand what a topological approach means, let's pretend that we have a rubbery object, like a circular rubber band. Topologists—that is, mathematicians who study topology—try to describe the properties of an object that stay the same if we stretch or shrink it or bend it, but without gluing things together, breaking the object, or creating any sharp points. From a topological viewpoint, because we can stretch the rubber band into an oval, we say that the circle and the oval are "equivalent." However, our rubber band is not equivalent to a segment of a string, because the rubber band has a hole in the middle but string does not. Remember that we're not allowed to glue the ends of the string together, and we're also not allowed to cut the rubber band.

By figuring out which shapes are equivalent to each other in this special way, we separate shapes into different groups. As an example, let's assign the letters in "Pokémon" to groups. The letters 'P' and 'o' belong to the same group, because we can compress the "tail" of the 'P' upwards and then stretch the hole into the shape of the letter 'o'. We illustrate this transformation in a short animation at https://drive.google.com/file/d/1SNrU0_usYBC7t5h_rHdkTvzHBC3-O5cf/view?usp=sharing. Consequently, 'P' and the two copies of 'o' make up one group of topologically equivalent letters. The 'k', 'm', and 'n' form a different group of topologically equivalent letters, because we can turn each of them into a disc by squeezing and bending them. Take a look again at our animation. The remaining letter, 'é', is an interesting one. Without its accent, we would be able to shrink the round tail of the 'e' into the left side of the semicircle at the top of the letter. We could then stretch that semicircle into the shape of the letter 'o', which places it into the same group as the two copies of 'o'. However, with the accent, 'é' has two separate pieces that we are not allowed to glue together, so it belongs to its own group. Look at our short animation once more.

Shapes that are in the same group have important things in common. Although the details of the shapes 'P' and 'o' differ, they both have one hole that we can't remove. By contrast, the letters 'k', 'm', and 'n' don't have any holes. If we take the uppercase letter `B', we see that it is not in either of these groups. It is, however, in the same group as the number '8', as both shapes have two holes. The number of pieces is also important, so the 'é' (with one hole and two pieces) is in a different group from all of the other letters that we've discussed. Try separating the letters in your name into topologically equivalent groups. Now let's make things even more interesting by looking at some Pokémon themselves. For each Pokémon in Fig. 1, count the numbers of pieces and holes. Are you able to group any of them together based on this information?

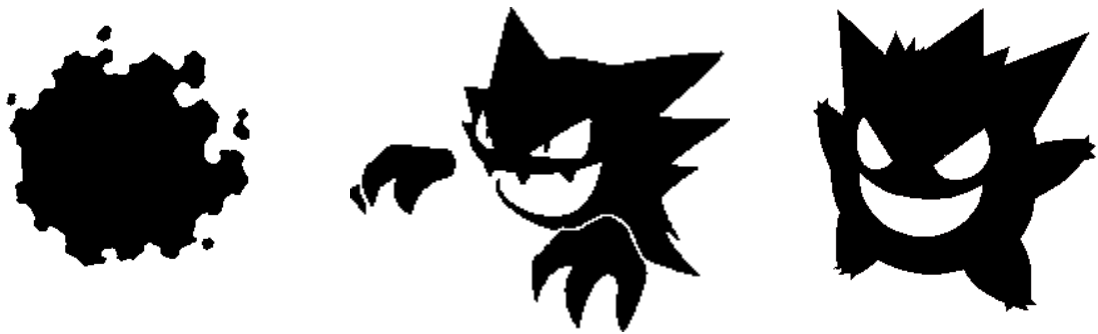

**Fig 1.** Pokémon have different shapes. [*Question*] Can you categorize Pokémon of different shapes based on their numbers of holes and numbers of pieces?

It can be challenging to study the topology of solid objects like those that we've been discussing so far, but now think about when you draw pictures in activities like "connect the dots." We see a bunch of dots—these dots are examples of *data* points—but we see enough of them that we often have a good idea of what shape we'll get when we connect the dots. We are good at determining shapes from just these dots, but is there a way to do this automatically? Take a look at the connect-the-dot activity for Pokémon in Fig. 2, and imagine what they look like when we connect the dots. This type of activity is often harder for a computer than for a human, but because we want to do this for many collections of dots — that is, when we have many thousands of games of connect the dots — mathematicians and other scientists seek good ways to do it automatically.

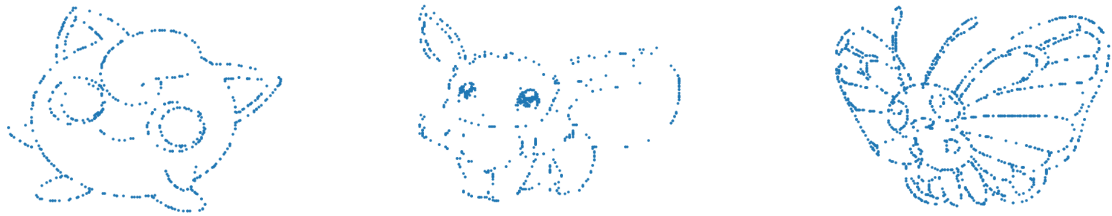

**Fig 2.** Connect the dots with Pokémon: (a) Jigglypuff, (b) Eevee, and (c) Butterfree.

As we'll see soon, topology can help us make sense of large amounts of data, and we can think of exploring the topology of a collection of data (called a *data set*) as a giant game of connecting the dots. In fact, *data* is a broad term that doesn't just come in the form of dots on a page. There are many different kinds of data, but we'll focus on data that have numbers attached to them, such as the locations of a bunch of dots on a page, the heights of the children at a school, or the number of words in each paragraph of this article.

## Discovering the shape of data

People think about topology and data together in an area of study that is called *topological data analysis* (TDA) [4,6,8]. In TDA, we try to describe the shape of data by first building a series of pictures by connecting the dots in the data. By connecting the dots in one of various ways (see [2] for a comparison of different ways in a study of voting data from the 2016 U.S. presidential election), we can use ideas from topology to study structures in the data. In particular, instead of connecting the dots by drawing lines from one dot to another like we're used to doing, let's connect the dots by increasing the size of the dots. As we make the dots larger and larger, the gaps between the dots become smaller and smaller, and eventually they overlap. Take a look at Jigglypuff in Fig. 3. The dots in Fig. 3(a) are tiny, but if we put a disc (that is, a filled-in circle) around each dot, we get Fig. 3(b). As we indicate in Table 1, the dotted Jigglypuff in Fig. 3(a) has 224 disconnected pieces but no holes. The drawing in Fig. 3(b) still has no holes, but now there are only 101 separate pieces.

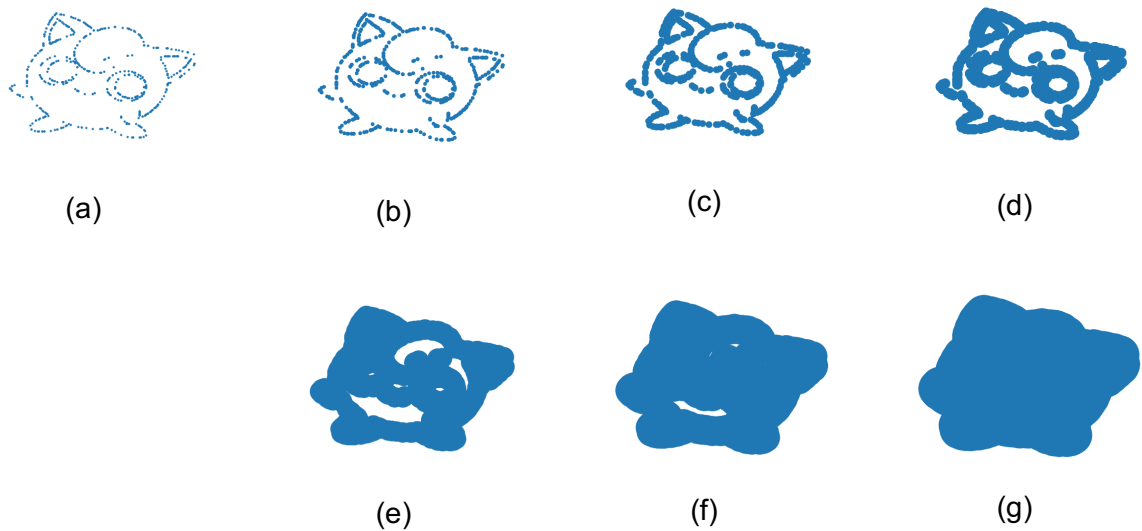

**Fig 3.** Increasing the size of dots in data. At first, Jigglypuff becomes easier to see as the dots get larger. Jigglypuff then gets harder to see.

But we aren't done yet, because we also need to figure out how large to make the discs, and we want to do this in a very thorough and careful way. Think of an extreme case where we make the discs really huge. We then get one very large object with no holes, as one can see for Jigglypuff in Fig. 3(g). However, we may see plenty of interesting things for different sizes of discs, ranging from very tiny ones (which may not be much larger than the dots themselves) to extremely large ones. The different versions of Jigglypuff in Fig. 3 have different numbers of pieces and different numbers of holes. By using mathematics and computation, we can consider many different sizes of discs—all sizes at once, in fact—and we obtain an object for each one. Each of the seven versions of Jigglypuff in Fig. 3 is one object. For each of these objects, we can count how many pieces it has and how many holes it has, just like in the pictures in Fig. 1.

Now let's go through the pictures of Jigglypuff in Fig. 3 in some more detail and summarize what we find in Table 1. When the discs are small, they don't touch each other, so there are no holes, as we indicate in the first row of Table 1. In Fig. 3(b), the discs have gotten a bit larger and some touch each other, so Jigglypuff now consists of fewer separate pieces. There are still no holes. In the third row of Table 1, which gives the count of pieces and holes in Fig. 3(c), even more discs are touching, so Jigglypuff consists of even fewer pieces. There are also now two holes, which appear to be Jigglypuff's left eye and right ear. As the discs get even larger, they combine with each other into one piece (see the second-last row of Table 1) with six holes. If we continue to let the discs grow, those holes fill in and we eventually obtain one large, hole-less blob; see Fig. 3(g) and the bottom row of Table 1. The information in Table 1 is one way of describing and summarizing what we observe from this range of disc sizes. That is, we are studying the structure of Jigglypuff across many size *scales*. Each version of Jigglypuff in Fig. 3 is at one scale, and by counting the number of pieces and number of holes at each scale, we can see the range of disc sizes over which Jigglypuff's features persist. People who study topological data analysis are very interested in how long different features persist in the data that they study.

| | Picture comes from | Number of pieces | Number of holes |
|---|---|---|---|
| 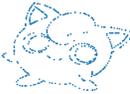 | Fig. 3(a) | 224 | 0 |
| 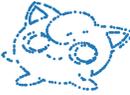 | Fig. 3(b) | 101 | 0 |
| 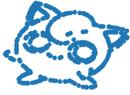 | Fig. 3(c) | 17 | 2 |
| 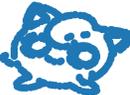 | Fig. 3(d) | 1 | 6 |
| 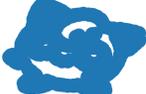 | Fig. 3(e) | 1 | 6 |
| 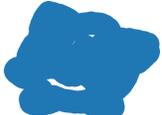 | Fig. 3(f) | 1 | 3 |
| 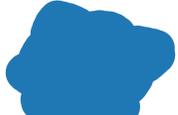 | Fig. 3(g) | 1 | 0 |

**Table 1.** We describe the six pictures of Jigglypuff in Fig. 3 with the numbers (224, 0), (101, 0), (17, 2), (1, 6), (1, 6), (1, 3), and (1, 0). In each pair, the first value indicates the number of

> distinct pieces and the second value indicates the number of holes.

## What can we learn from topological data analysis?

Topological data analysis — a mathematical way of building on ideas like "connecting the dots" — can tell us a lot about the world and many things in it. It allows us to explore complex data from a huge variety of topics in social science, biology, astronomy, and more [4]. Although mathematics and computation give many ways for people to study data, TDA is designed specifically to help understand shape and scale.

We can use TDA to help us understand the universe. Planets like Earth are part of solar systems, which in turn are part of galaxies, which occur in clusters. If we look into a telescope and zoom in on a solar system, the planets seem to be very far apart. But if we zoom out to look at an entire galaxy, each solar system may appear as just a dot and structures within the solar system seem to merge together. If we zoom out further, each galaxy may appear as just a dot. To study the structure of the universe at these different scales, scientists have used TDA to count the number of pieces and holes in a data set of star positions [1].

Moving back down to earth, scientists have used TDA to examine the patterns of veins in leaves [5]. They studied the structure of more than 100 leaves and used topology to find different characteristics (kind of like human fingerprints) of patterns in the leaves. These fingerprints, which give different information from other ways of studying these patterns, improve the ability of scientists to identify leaves from small fragments of them and suggest a possible mechanism for how leaves grow. A similar approach is useful for studying the structure of fungi, blood vessels, and other things with branches and loops.

People also use TDA to describe activity patterns of people, animals, and computer simulations of people and animals. In [2], for example, two of us studied the structure of voting in different areas of California. We used TDA to detect areas of the state where people voted differently in the 2016 presidential election from those in neighboring areas. Animals other than people also produce interesting patterns. Animal groups, such as schools of fish and flocks of birds, have many individuals and can form beautiful structures. TDA can help scientists explore and understand these complex patterns, like in recent research that uses simulations of flocks [9].

Topological data analysis is an increasingly popular approach for studying many problems, ranging from "connecting the dots" in pictures of Pokémon to the structure of the universe [1], patterns in nature [5], who people vote for in elections [2], and much more. TDA is a fascinating and important area of mathematics for helping to make sense of complex data [4,6,8].

# Glossary

**Data**: Information that is collected. A collection of data is sometimes called a *data set*.

**Scale**: A characteristic size of an object, such as the radius of a disc or the length of a side of a square.

**Topologically equivalent**: A term that describes two objects that can be turned into each other by stretching, shrinking, bending, or warping them (but not gluing or tearing them).

**Topology**: A branch of mathematics that people use to study the shapes of objects. People who study topology are sometimes called *topologists*.

# Answer Key

**Question.** (a) Gastly has 5 pieces and 0 holes. (b) Haunter has 4 pieces and 2 holes. (c) Gengar has 1 piece and 3 holes. To be in the same group, two objects (which, in this case, are Pokémon) should have both the same number of pieces and the same number of holes. Because Gastly, Haunter, and Gengar have different numbers of pieces and different numbers of holes, they are in three different groups.

## Conflict-of-Interest Statement

The authors declare that the research was conducted in the absence of any commercial or financial relationships that could be construed as a potential conflict of interest.

## Acknowledgements

We are grateful to our young readers — Charlotte Amann-Sulzmann, Simon Cafiero, Addison Cart, Nia Chiou, Valerie K. Eng, Linnea Keiser-Clark, Coralea Lash-St. John, Adele Low, Maple Leung, Nora Stricker, Kate Van Hooser, and one anonymous person — for their many helpful comments. We also thank their parents, teachers, and friends — Clayton Cafiero, Lyndie Chiou, Puck Rombach, and Steve Van Hooser — for putting us in touch with them and soliciting their feedback. We also thank Norman Redington and our reviewers for helpful comments. MAP, MF, and YHK acknowledge support from the National Science Foundation (grant number 1922952) through the Algorithms for Threat Detection (ATD) program. CMT acknowledges support from the National Science Foundation (grant number 1813752) through the Division of Mathematical Sciences.

## Author Biographies

**Michelle Feng** is a postdoctoral researcher in Computing + Mathematical Sciences at Caltech. Her research focuses on studying the shape of political and social data. She is interested in social justice and understanding the construction of social histories. In addition to her mathematical interests, Michelle is a fan of video games and animation. She's a huge Pokémon fan (some might call her a Pokéfanatic), and she loves ghost and fire types, but especially Vulpix, Phantump, and Fennekin.

**Abigail Hickok** is a second-year graduate student in Mathematics at UCLA. She grew up in New Jersey, but is learning to love the west coast, too. Her research focus is topological data analysis. When she's not doing math, she likes to spend time outdoors by going rock climbing and hiking, and she also enjoys playing piano. Her favorite Pokémon is Charmander.

**Yacoub H. Kureh** was born and raised in Orange County, California and recently earned his Ph.D. in Mathematics at UCLA. His research interests are in the areas of network science and data science. He focuses on modeling opinions spreading among individuals who share social connections. He is a first-generation college and graduate student. He is passionate about education and supporting the right to education. Yacoub fondly remembers playing the original Pokémon Red and Pokémon Blue games with his grade-school friends. His favorite Pokémon are Ditto and Eevee.

**Mason A. Porter** is a professor in the Department of Mathematics at UCLA. He was born in Los Angeles, California, and he is excited to be a professor in his hometown. In addition to studying networks and other topics in mathematics and its applications, Mason is a big fan of games of all kinds, fantasy, baseball, the 1980s, and other delightful things. Mason used to be a professor at University of Oxford, where he did actually wear robes on occasion (like in the *Harry Potter* series). When not isolating himself from the rest of humanity, Mason occasionally likes to spend time with other people. His Pokémon of choice is Jigglypuff (obviously).

**Chad M. Topaz** is Professor of Mathematics at Williams College in Massachusetts. He is fascinated by patterns in nature and has used mathematics to study the patterns that form in animal swarms, chemical reactions, waves in fluids, and more. Chad was born in Chicago and has also lived in Boston, Raleigh-Durham, San Francisco, Los Angeles, and Minneapolis before landing in Williamstown, home to seven thousand people and one stoplight (but the stoplight is pretty unnecessary). After writing this paper, Chad feels bad that he never tagged along with his husband and daughter on their Pokémon Go outings.